# The Square of the Dirichlet-to-Neumann map equals minus Laplacian

D V Ingerman

**Abstract.** The Dirichlet-to-Neumann maps connect boundary values of harmonic functions. It is an amazing fact that the square of the non-local Dirichlet-to-Neumann map for the uniform conductivity 1 on the unit disc equals minus the local(!) Laplace operator on the boundary circle. To establish a new connection between discrete and continuous Dirichlet-to-Neumann maps and for the approximations I construct a finite and an infinite graphs which Dirichlet-to-Neumann map have the same property: Λ²(1) = - δ²/δθ². The construction gives a new continued fraction identity. It is interesting to consider the geometric and probabilistic (trajectories of the random walk) consequences of this localizing identity unifying discrete and continuous equations for potentials.

**1.Introduction**
The continuous Dirichlet-to-Neumann map is a map from boundary potential to boundary current of a body interior potential satisfying Laplace equation.
$$(\gamma u_x)_x + (\gamma u_y)_y = 0.$$

The discrete Dirichlet-to-Neumann map of a graph is the linear map from boundary potential to boundary current for a potential on vertices of the graph satisfying Kirchoff's Law.
$$\sum_j \gamma_{ij}(u_i - u_j) = 0.$$

These Dirichlet-to-Neumann maps are the data for inverse boundary problems arising in the problems requiring non-intrusive testing in medicine and oil and gas production industry. There are discrete (on graphs) and continuous (on manifolds) models for the inverse problems. There are several connections found between discrete and continuous Dirichlet-to-Neumann maps in the two dimensional case:
- restrictions of the kernel of the continuous Dirichlet-to-Neumann maps are the matrices that are Dirichlet-to-Neumann maps of circular planar graphs, see [4] and [3],
- spectrum interpolation functions of the continuous Dirichlet-to-Neumann maps in the radially symmetric layered case are the limits of the continued fractioins representing Dirichlet-to-Neumann maps of the graphs (discrete discs), see [2].

Professor Uhlmann has noticed a fundamental property of the **non-local** Dirichlet-to-Neumann map for the uniform conductivity *1* on the unit disc: its square equals **local** operator, the minus Laplacian on the boundary circle.
$$\Lambda^2(1) = -\frac{\partial^2}{\partial \theta^2}.$$

He then has asked the question if the similar property is true in the discrete case. In this paper I formulate a discrete analog of the fundamental property for graphs Γ

$$\Lambda^2(\Gamma) = \begin{bmatrix} 2 & -1 & 0 & 0 & \ldots & 0 & -1 \\ -1 & 2 & -1 & 0 & \ldots & 0 & 0 \\ 0 & -1 & 2 & -1 & \ldots & 0 & 0 \\ \ldots & \ldots & \ldots & \ldots & \ldots & \ldots & \ldots \\ 0 & 0 & 0 & \ldots & \ldots & -1 & 0 \\ 0 & 0 & 0 & \ldots & -1 & 2 & -1 \\ -1 & 0 & 0 & \ldots & 0 & -1 & 2 \end{bmatrix}$$

and give two explicit constructions of a finite and infinite graphs with the discrete fundamental property. This establishes a new connection between discrete and continuous Dirichlet-to-Neumann maps and implies many determinant and eigenvalue inequalities for blocks of the **full** matrix that is a square root of the **sparse** tridiagonal matrix, see [2] and [3]. The construction gives a new continued fraction identity:

$$\frac{w^k + w^{-k}}{w^k - w^{-k}}(w^l - w^{-l}) + \cfrac{1}{\cfrac{w^{k-1} + w^{1-k}}{w^{k-1} - w^{1-k}}(w^l - w^{-l}) + \ldots + \cfrac{1}{\cfrac{w^2 + w^{-2}}{w^2 - w^{-2}}(w^l - w^{-l}) + \cfrac{1}{\cfrac{w + w^{-1}}{w - w^{-1}}(w^l - w^{-l})}}} = 1.$$

where

$$w = e^{\frac{i\pi}{2k+1}}.$$

It is interesting to consider the geometric and probabilistic (trajectories of the random walk [5]) consequences of this localizing identity unifying discrete and continuous equations for potentials.

## 2. Continuous Case

The Dirichlet-to-Neumann map $\Lambda(\gamma)$ for a positive function $\gamma$ (conductivity) on a unit disc is the linear map

$$\Lambda(\gamma): f \to g$$

where $u$ (potential) is the solution of the following Laplace equation

$$\partial(\gamma \partial u / \partial x)/\partial x + \partial(\gamma \partial u / \partial y)/\partial y = 0$$
$$f(\theta) = u(1, \theta)$$

is the restriction of $u$ on the boundary circle and
$$g(\theta) = \gamma \partial u(r, \theta)/\partial r \,|\, r = 1$$
is the boundary current.
For $\gamma = 1$

$$u(z) = z^l$$
$$u(r, \theta) = r^l \sin(l\theta), u(r, \theta) = r^l \cos(l\theta)$$
$$\Lambda(1)(\sin(l\theta) = l\sin(l\theta), \Lambda(1)(\cos(l\theta)) = l\cos(l\theta)$$

and therefore the square of the **non-local** Dirichlet-to-Neumann map for the uniform conductivity $1$ on the unit disc equals minus the **local** Laplace operator on the boundary circle

$$\Lambda^2(1)(\sin(l\theta)) = l^2 \sin(l\theta)$$
$$\Lambda^2(1)(\cos(l\theta)) = l^2 \cos(l\theta)$$
$$\Downarrow$$
$$\Lambda^2(1) = -\frac{\partial^2}{\partial \theta^2}.$$

## 3. Discrete Case: Finite Graph

By analogy with the continous case of a body the discrete Dirichlet-to-Neumann map of a graph is the linear map from boundary potential to boundary current for a potential on vertices of the graph satisfying Kirchoff's Law. For a definition of the Dirichlet-to-Neumann map as a Schur complement of the adjacency/Laplace matrix please see [3]. The discrete analog of the local Laplace operator on the boundary circle of the unit disc is the sparse tridiagonal matrix L or *[- ∂²/δθ ²]* in the notation of [2]:

$$L = \left[-\frac{\partial^2}{\partial \theta^2}\right] = \begin{bmatrix} 2 & -1 & 0 & 0 & \ldots & 0 & -1 \\ -1 & 2 & -1 & 0 & \ldots & 0 & 0 \\ 0 & -1 & 2 & -1 & \ldots & 0 & 0 \\ \ldots & \ldots & \ldots & \ldots & \ldots & \ldots & \ldots \\ 0 & 0 & 0 & \ldots & \ldots & -1 & 0 \\ 0 & 0 & 0 & \ldots & -1 & 2 & -1 \\ -1 & 0 & 0 & \ldots & 0 & -1 & 2 \end{bmatrix}$$

$L = \{L(i, i) = 2, L(i, i+1) = L(i, i-1) = -1, \text{circularly}\} = [-\partial^2/\delta\theta^2]$.

The existence of the finite discrete graph $\Gamma$ with the desired fundamental property

$$\Lambda^2(\Gamma) = L = \left[-\frac{\partial^2}{\partial \theta^2}\right].$$

was proved in [2]. In this paper I give an explicit formula for the discrete parameters (conductivities) on the graph.

It follows from [2] that to construct a finite network (discrete disc) $\Gamma$

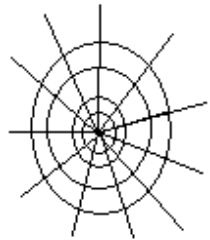

**Figure 1.** Radially symmetric network $\Gamma$. We found explicit formula for conductivities on $\Gamma$ that give the fundamental property in the title.

with the desired property

$$\Lambda^2(\Gamma) = \left[-\frac{\partial^2}{\partial \theta^2}\right].$$

it is enough to find the coefficients $c_i$ (layer conductivities) in the continued fraction

$$\beta(\lambda) = c_k \lambda + \cfrac{1}{c_{k-1}\lambda + \cfrac{1}{c_{k-2}\lambda + \ldots + \cfrac{1}{c_2\lambda + \cfrac{1}{c_1\lambda}}}} = 1$$

for the square root of every non-zero eigenvalue $l$ of the *2k+1* by *2k+1* discrete minus Laplacian matrix $L$

$$L = \left[-\frac{\partial^2}{\partial \theta^2}\right] = \begin{bmatrix} 2 & -1 & 0 & 0 & \ldots & 0 & -1 \\ -1 & 2 & -1 & 0 & \ldots & 0 & 0 \\ 0 & -1 & 2 & -1 & \ldots & 0 & 0 \\ \ldots & \ldots & \ldots & \ldots & \ldots & \ldots & \ldots \\ 0 & 0 & 0 & \ldots & \ldots & -1 & 0 \\ 0 & 0 & 0 & \ldots & -1 & 2 & -1 \\ -1 & 0 & 0 & \ldots & 0 & -1 & 2 \end{bmatrix}$$

The direct calculation shows that for

$$l = 1, 2, \ldots, k$$

$$\lambda = w^l - w^{-l}$$

where $w$ is the *(2k+1)*st root of unity

$$w^{2k+1} = -1$$

such that

$$w \neq -1.$$

**Conjecture 3.1**

$$\frac{w^k + w^{-k}}{w^k - w^{-k}}(w^l - w^{-l}) + \cfrac{1}{\cfrac{w^{k-1} + w^{1-k}}{w^{k-1} - w^{1-k}}(w^l - w^{-l}) + \ldots + \cfrac{1}{\cfrac{w^2 + w^{-2}}{w^2 - w^{-2}}(w^l - w^{-l}) + \cfrac{1}{\cfrac{w + w^{-1}}{w - w^{-1}}(w^l - w^{-l})}}} = 1.$$

for

$$l = 1, 2, \ldots, k.$$

where

$$w = e^{\frac{i\pi}{2k+1}}.$$

That is

$$c_l = \frac{w^l + w^{-l}}{w^l - w^{-l}}.$$

I have guessed Conjecture *3.1* on a calculator and a computer and checked that it is true for up to *41* floor continued fraction on *MATLAB*. For small $k$ it has a simple geometric meaning involving relationships between lengths of diagonals, heights and medians of symmetric polygons.

Please notice that even for *k = 1* the conjectured identity (expression $\beta(\lambda) = 1$) involves *5* fundamental operations and *4* fundamental constants of mathematics: addition (+), subtraction (-), multiplication (*), division (/), power (^), *e, π, i, 1*.

Zeros and poles of the continued fraction $\beta(\lambda)$ and polynomials $p(\lambda)$ and $q(\lambda)$, where $\beta(\lambda) = p(\lambda)/q(\lambda)$, have a similar to zeros of the Riemann $\zeta$-function property, they lie on a single line [2].

## 4. Discrete Case: Infinite Graph

Let *S* be a set of *n* points on a circle. The vertices of the constructed infinite graph $\Gamma$ are the set *SxN*, where *N = 1, 2, 3...* with all edges of conductivity *1* connecting vertex *i* of the layer *Sx{j}* to two vertices of the next layer *i* and *i + 1* of *Sx{j+1}* circularly. In other words $\Gamma$ is *ZxZ* cut along a diagonal *x = y* and rolled into a semi-infinite cylinder. This graph is self-dual and self-similar.

Here are pictures of $\Gamma$ for *n=3,5*.

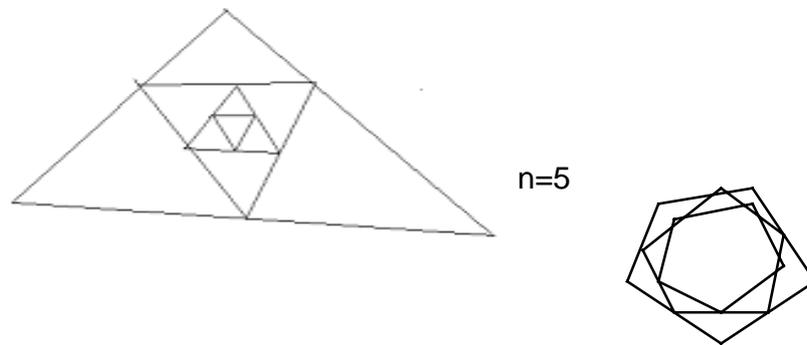

**Figure 2.** Infinite graphs with the desired fundamental property in the title with 3 and 5 boundary vertices.

**Theorem 4.1**

$$\Lambda^2(\Gamma) = L = \left[-\frac{\partial^2}{\partial \theta^2}\right]$$

or

$$\Lambda^2(\Gamma) = \begin{bmatrix} 2 & -1 & 0 & 0 & \ldots & 0 & -1 \\ -1 & 2 & -1 & 0 & \ldots & 0 & 0 \\ 0 & -1 & 2 & -1 & \ldots & 0 & 0 \\ \ldots & \ldots & \ldots & \ldots & \ldots & \ldots & \ldots \\ 0 & 0 & 0 & \ldots & \ldots & -1 & 0 \\ 0 & 0 & 0 & \ldots & -1 & 2 & -1 \\ -1 & 0 & 0 & \ldots & 0 & -1 & 2 \end{bmatrix}$$

*Proof:*
The Dirichlet-to-Neumann maps of a graph on the figure *2* with an added layer is the same as of the original graph since

$$N + \{0\} = N$$

According to [3] the Dirichlet-to-Neumann map $\Lambda$ of $\Gamma$ is the Schur complement of the Laplace matrix, as the fixed point:

$$\Lambda = \begin{bmatrix} 2I & B^T \\ B & \Lambda + 2I \end{bmatrix} / [\Lambda + 2I]$$

therefore

$$\Lambda = 2I - B^T(\Lambda + 2I)^{-1}B$$

where $B$ is the adjacency matrix

$$B = \begin{bmatrix} -1 & -1 & 0 & 0 & \ldots & 0 & 0 \\ 0 & -1 & -1 & 0 & \ldots & \ldots & 0 \\ 0 & 0 & -1 & -1 & \ldots & \ldots & 0 \\ \ldots & \ldots & \ldots & \ldots & \ldots & \ldots & \ldots \\ 0 & 0 & \ldots & \ldots & -1 & -1 & 0 \\ 0 & 0 & \ldots & \ldots & \ldots & -1 & -1 \\ -1 & 0 & 0 & \ldots & \ldots & 0 & -1 \end{bmatrix}$$

$B = \{B(i, i) = B(i, i+1) = -1, \text{circularly}\}$.

Therefore

$$(2I - \Lambda)(2I + \Lambda) = BB^T$$

and

$$\Lambda^2 = 4I - BB^T = 4I - |L| = L = [-\partial^2/\partial\theta^2].$$

According to [3] the constructed finite and infinite graphs are $Y$-$\Delta$ transformation equivalent since they are both circular planar and have equal Dirichlet-to-Neumann maps. It is interesting to consider applications of this fact to transformations of finite and infinite graphs, networks, knots and Boolean circuits.